\documentclass{amsart}
\usepackage[dvips]{graphicx}
\newcommand{\re}{{I\!\!R}}
\newcommand{\ren}{\re^N}

\newcommand{\dyle}{\displaystyle}
\newcommand{\ene}{{I\!\!N}}

\newcommand{\irn}{\int_{\re^N}}
\newcommand{\io}{\int_{\O}}

\newcommand{\limit}{\lim\limits}

\newcommand{\dob}{{\mathcal D}^{1,2}}

\newcommand{\db}{\rightharpoonup}

\renewcommand{\a }{\alpha }

\renewcommand{\d }{\delta }
\newcommand{\D }{\Delta }
\newcommand{\e }{\varepsilon }
\newcommand{\nab }{\nabla }

\newcommand{\inr }{\int _{\ren} }

\newcommand{\g }{\gamma}
\newcommand{\G }{\Gamma }
\renewcommand{\l }{\lambda }
\renewcommand{\L }{\Lambda }

\newcommand{\n }{\nabla }

\renewcommand{\O }{\Omega }

\newenvironment{pf}{\noindent{\sc Proof}.\enspace}{\rule{2mm}{2mm}\medskip}
\newtheorem{Theorem}{Theorem}[section]
\newtheorem{Corollary}[Theorem]{Corollary}
\newtheorem{Definition}[Theorem]{Definition}
\newtheorem{Lemma}[Theorem]{Lemma}

\newtheorem{remark}[Theorem]{Remark}
\begin{document}

\title[An  equation involving Hardy inequality
and critical Sobolev exponent]{Existence and multiplicity for
perturbations of an  equation involving Hardy inequality and
critical Sobolev exponent in the whole $\ren$.}

\author[B. Abdellaoui, V. Felli \& I. Peral]{B. Abdellaoui, V. Felli,  \& I.
Peral}

\address{{\bf B. Abdellaoui\& I. Peral}, Departamento de Matem\'aticas, U. Autonoma de Madrid,
\hfill\break\indent 28049 Madrid, Spain.} \email{{\tt
boumediene.abdellaoui@uam.es, ireneo.peral@uam.es}}

\address{{\bf V. Felli}, S.I.S.S.A., Via Beirut 2-4, 34014 Trieste, Italy}
\email{{\tt felli@sissa.it} }

\thanks{First and third authors supported by project BFM2001-0183, M.C.Y.T.
Spain.\, Second author supported by Italy MIUR, national project
``Variational Methods and Nonlinear Differential Equations''}


\begin{abstract}
\noindent In order to obtain solutions to problem
$$
\left\{
\begin{array}{c}
-\D u=\dfrac{A+h(x)} {|x|^2}u+k(x)u^{2^*-1},\,\,x\in \ren,
\\ u>0 \hbox{  in  }\ren, \mbox{  and  }u\in \dob(\ren),
\end{array}
\right. $$ $h$ and $k$ must be chosen  taking into account not
only the size of some norm but the shape. Moreover, if $h(x)\equiv
0$, to reach multiplicity of solution, some hypotheses about the
local behaviour of $k$ close to the points of maximum are needed.
\end{abstract}
\maketitle
\date{}
\section{Introduction}\label{s:0}
In this paper we will consider the following class of problems
\begin{equation}\label{eq:poi0}
\left\{
\begin{array}{c}
-\D u=\Big(\dfrac{A +h(x)}{|x|^2}\Big)u+k(x)u^{2^*-1},\,\,x\in \ren,
\\ u>0 \hbox{  in  }\ren, \mbox{  and  }u\in \dob(\ren),
\end{array}
\right.
\end{equation}
where $N\ge 3$,  $2^*=\frac{2N}{N-2}$ and $h,\, k$ are continuous
bounded functions, for which we will state appropriate
complementary hypotheses. Here  $\dob(\ren)$ denotes the closure
space of ${\mathcal C}_0^\infty(\ren)$ with respect to the norm
$$
||u||_{\dob(\ren)}:=\bigg(\irn |\n u|^2dx\bigg)^{1/2}.
$$
By the
Sobolev inequality we can see that $\dob(\ren)$ is the class of
functions in $L^{2^*}(\ren)$ the distributional gradient of which
 satisfies $\big(\irn |\n
u|^2dx\big)^{1/2}<\infty$.

\

\noindent For $h\equiv 0, k\equiv 1$ the problem is studied  by S.
Terracini in \cite{S}. In \cite{FC}  the existence of a positive
solution  is proved in the case $h=0$  by using the perturbative
method by Ambrosetti-Badiale in \cite{AB},  even for a more general class of
differential operators related to the Caffarelli-Kohn-Nirenberg
inequalities that  contains our operator. By the perturbative
nature of the method, the solutions found in \cite{FC}  are close
to some radial solutions to  the unperturbed problem. On the other
hand, in \cite{SM} Smets obtains the existence of a positive
solution for  problem  \eqref{eq:poi0}   with $h=0$,   $k$
bounded, $k(0)=\limit\limits_{|x|\to\infty}k(x)$ and dimension
$N=4$.

\noindent In this paper we study the existence of positive solutions
in the case in which either $h\equiv 0$ and $k\not\equiv1$ or $k\equiv
1$ and $ h\not\equiv 0$ satisfying suitable assumptions. Our results hold
in any dimension and are proved using the {\it
  concentration-compactness} arguments by P.L. Lions.

\noindent It is known that the general problem has an obstruction
provided by a Pohozaev type identity that shows us the
particularity of this problem, that is: {\it the existence of a
positive solution depends not only on the size of the functions
$h$ and $k$ but also on their shape}. More precisely, assume that
$u$ is a variational solution to our equation with $h,k\in
\mathcal{C}^1$. Multiplying the equation by $\langle x,\nabla
u\rangle$ and with a convenient argument of approximation we get
that necessarily $$\frac {\l}2\int \langle \nabla
h(x),x\rangle\frac{u^2}{|x|^2}dx+
 \frac {1}{2^*}\int \langle \nabla k(x),x\rangle|u|^{2^*}dx=0.$$
This  behaviour makes the problem more interesting to be analyzed.
The existence part of the paper is mainly based on the {\it
concentration-compactness arguments} by P.L. Lions (see \cite{PL1}
and \cite{PL2}) and involves some qualitative properties of the
coefficients that avoids the Pohozaev type obstruction. We also
obtain  multiplicity of positive solutions by using variational
and topological arguments.

The organization of the paper is as follows. Section 2
is devoted to the study of nonexistence and existence for $k\equiv
1$ and $h$ satisfying suitable conditions. As pointed out above,
we mainly use the {\it concentration-compactness} principle by
P.L. Lions. The main result in this part is Theorem \ref{th:ex}.
 Section 3 deals with the existence and multiplicity
results for the case in which $h\equiv0$ and $k$ satisfies some
convenient conditions. In this part of the paper we will use
techniques that previously had been introduced to study related
problems by Tarantello in \cite{TA} and refined by Cao- Chabrowsky
in \cite{CC1} (see also the references therein).We use this
approach in the case that the function $k$ achieves its maximum at
a finite number of points. The main result in Section 3 is
Theorem \ref{th:mutlip}.  In Section~4 we  study a more general
class of functions $k$, i.e. we treat the case in which $k$ can reach its
maximum at infinitely many points, but having only accumulation
points at finite distance to the origin. To analyze this case we
use  the  {\it Lusternik-Schnirelman category}. This point of view
is inspired by the study of multiplicity of positive solutions to
subcritical problems done by R. Musina in \cite{Mu}. After several
technical lemmas the main result contained in Section 4 is Theorem
\ref{th:ultimo}.

In a forthcoming paper we will discuss the case of critical
equations related to the so called Caffarelli-Kohn-Nirenberg
inequalities.

\

{\bf Acknowledgment.-} We want to thank Professor A. Ambrosetti
for his encouragement and for many helpful  suggestions. Part of this
 work was carried out while the second
author was visiting Universidad Aut\'{o}noma of Madrid; she wishes
to express her gratitude
to Departamento de Matem\'{a}ticas of Universidad Aut\'{o}noma
for its warm hospitality.

\section{Perturbation in the linear term}\label{s:1}

We will study perturbations of a  class of elliptic
equations in $\ren$ related to a Hardy inequality interacting with
a nonlinear term  involving the critical Sobolev exponent .
Precisely we will consider the following problem
\begin{equation}\label{eq:poi1}
\left\{
\begin{array}{c}
-\D u=\dfrac{A +h(x)}{|x|^2}u+u^{2^*-1},\,\,x\in \ren,
\\ u>0 \hbox{  in  }\ren, \mbox{  and  }u\in \dob(\ren),
\end{array}
\right.
\end{equation}
where $N\ge 3$ and $2^*=\frac{2N}{N-2}$. Hypotheses on $h$ will be
given below. To be precise  we recall the Hardy inequality.
\begin{Lemma}\label{lm:hardy}{\it (Hardy inequality)}
Assume that $u \in \dob(\ren)$, then  $\dfrac{u}{|x|}\in L^2(\ren
)$ and $$\int_{\ren } \dfrac{|u|^2}{|x|^2}dx \le
C_{N}\int_{\ren}|\nabla u |^2dx,$$ where
$C_{N}=\big(\frac{2}{N-2}\big)^2$ is optimal and not attained.
\end{Lemma}
Hereafter we will call $\L_N:=C_{N}^{-1}=\frac{(N-2)^2}{4}$.
See for instance \cite{GP} for a proof.

\noindent The case $h=0$ of \eqref{eq:poi1} has been studied by S.
Terracini in \cite{S}; she shows, in particular, that
\begin{enumerate}
\item if $A \ge \L_N$, then problem (\ref{eq:poi1}) has no positive
solution in $\mathcal{ D}'(\ren)$;
\item if $A\in (0,\L_N)$ then problem (\ref{eq:poi1}) has the
one-dimensional $\mathcal{C}^2$ manifold of positive solutions
\begin{equation}\label{eq:ZA}
Z_A=\bigg\{ w_\mu\,\,|\,w_{\mu }(x)=\mu
^{-\frac{N-2}{2}}w^{(A)}\Big(\dfrac{x}{\mu }\Big),\,\,
\mu>0\bigg\},
\end{equation}
where
\begin{equation}\label{eq:zzz}
w^{(A)}(x)=\frac{\left(N(N-2)\nu
_A^2\right)^{\frac{N-2}{4}}}{\left(|x|^{1-\nu _A}(1+|x|^{2\nu
_A})\right)^\frac{N-2}{2}}, \mbox{ and }\nu
_A=\left(1-\frac{A}{\L_N}\right)^{\frac12}.
\end{equation}
Moreover, if we set $Q_A(u)=\dyle\inr |\nab u|^2dx -A\inr
\frac{u^2}{|x|^2}dx$, then we obtain that
\begin{equation}\label{eq:minimiz}
\overline{S}\equiv\inf_{u\in \dob(\ren)\backslash
\{0\}}\dfrac{Q_A(u)}{\|u\|_{2^*}^2}=\dfrac{Q_A(w_{\mu})}{\|w_{\mu}\|_{2^*}^2}
=S\bigg(1-\frac{A}{\L_N}\bigg)^{\frac{N-1}{N}},
\end{equation}
where $S$ is the best constant in the Sobolev inequality. Notice
that $\overline{S}$ is attained
 exactly in the family $w_{\mu}$ defined in \eqref{eq:ZA}.
\end{enumerate}
\subsection{Nonexistence results}\label{ss11}
\noindent We begin by proving some nonexistence results that show
the fact that in this kind of problems both the size and the shape
of the perturbation are important.
Define
\begin{equation}\label{eq:energia} Q(u)=\irn|\nab u|^2dx-\inr\bigg(\frac{A+
h(x)}{|x|^2}\bigg)u^2dx,
\end{equation}
${\mathcal K}=\big\{u\in \dob(\ren)\,|\,\inr |u|^{2^*} dx=1\big\}$, and consider $
I_1=\inf_{u\in{\mathcal K}}Q(u).$
\begin{Lemma}\label{lm:ne} Problem (\ref{eq:poi1}) has
no positive solution in the following cases:
\begin{enumerate}
\item[(a)] If $A+
h(x)\ge 0$ in some ball $B_{\d }(0)$ and $I_1< 0$.
\item[(b)] If
$h$ is a differentiable function such that $\langle h'(x),x
\rangle$ has a fixed sign.
\end{enumerate}
\end{Lemma}
\begin{pf}
We begin by proving nonexistence under hypothesis $(a)$. Suppose
that $I_1<0$, and let $u$ be a positive solution to
(\ref{eq:poi1}). By classical regularity results for elliptic
equations we obtain that $u\in
\mathcal{C}^{\infty}(\ren\backslash\{0\})$. On the other hand,
since $A+ h(x)\ge 0$ in $B_{\d}(0)$, we obtain that $ -\D u\ge 0
\mbox{ in }\mathcal{D}'(B_{\d}(0))$. Therefore, since $u\ge 0$ and
$u\ne 0$, by the strong maximum principle we obtain that $u(x)\ge
c>0$ in some ball $B_{\eta }(0)\subset \subset
B_{\d}(0)$.
\\ Let $\phi _n\in \mathcal{C}^{\infty}_0(\ren)$, $\phi
_n\ge 0$, $||\phi_n||_{2*}=1$,  be a minimizing sequence of $I_1$.
By using $\dfrac{\phi_n^2}{u}$ as a test function in equation
$(\ref{eq:poi1})$ we obtain $$ \irn\nab \Big(\frac{\phi
_n^2}{u}\Big)\n u=\inr\frac{A+ h(x)}{|x|^2}\phi^2_n +\irn \phi
_n^2u^{2^*-2}. $$ A direct computation gives
\begin{equation}\label{eq:m30}
2\irn\frac{\phi _n}{u}\nab \phi _n\n udx - \irn\frac{\phi
_n^2}{u^2}|\nab u|^2dx = \inr\frac{A+ h(x)}{|x|^2}\phi^2_n+\irn
\phi _n^2u^{2^*-2},
\end{equation}
and since $$ 2\frac{\phi _n}{u}\nab \phi _n\n u-\frac{\phi
_n^2}{u^2}|\nab u|^2\le |\n \phi_n|^2,$$ we conclude that $$
\irn|\n \phi _n|^2dx \ge \inr\frac{A+ h(x)}{|x|^2}\phi^2_n+\irn
\phi _n^2u^{2^*-2}. $$ On the other hand, $I_1<0$ implies that we
can find an integer $n_0$ such that if $n\ge n_0$,
 $$ \irn|\nab \phi _n|^2-\inr\frac{A+ h(x)}{|x|^2}\phi^2_n <0. $$
As a consequence $\irn \phi_n^2u^{2^*-2}<0$, for $n\ge n_0$, which
 contradicts
the hypothesis $u>0$.

\noindent Let us now prove $(b)$. By using the Pohozaev  multiplier
$\langle x,\nabla u\rangle$, we obtain that if $u$ is a positive
solution to (\ref{eq:poi1}), then $$\irn \frac{\langle h'(x),x\rangle }{|x|^2}u^2dx=0,$$ which is not
possible if $\langle h'(x),x\rangle$ has a fixed sign and $u\not\equiv 0$.
\end{pf}

\begin{Corollary}
Assume either
\begin{itemize}
\item[$i)$] $A>\L_N$ and  $h\ge
0$,
or
\item[$ii)$]$A > \L_N$ and $1\le
\frac{4A}{(N-2)^2\|h\|_{\infty}}$,
\end{itemize}
then problem  (\ref{eq:poi1}), has no positive solution.
\end{Corollary}

\subsection{The local Palais-Smale condition: existence results}
\noindent To prove the existence results we will use a
variational approach for the associated functional
\begin{equation}\label{eq:energy} J(u)=
\frac12\irn |\n u|^2dx-\frac12\irn
\frac{A+h(x)}{|x|^2}u^2dx-\frac{1}{2^*}\irn |u|^{2^*}dx.
\end{equation}
We suppose that $h$ verifies the following hypotheses
\begin{enumerate}
\item[(h0)] $A+h(0)>0$.
\item[(h1)] $ h\in C(\ren)\cap L^{\infty}(\ren)$.
\item[(h2)] For some $c_0>0$, $A+||h||_{\infty}\le \L_N-c_0$.
\end{enumerate}
Critical points of $J$ in $\dob(\ren)$ are solutions
to equation (\ref{eq:poi1}). We begin by proving  a local Palais-Smale
condition for $J$. Precisely we prove the following Theorem.
\begin{Theorem}\label{th:ps}
Suppose that $(h0),\,(h1),\ (h2)$ hold and denote $h(\infty)\equiv
\limsup\limits_{|x|\to \infty}h(x)$.\\ Let $\{u_n\}\subset \dob(\ren)$
be a Palais-Smale sequence for
$J$, namely $$ J(u_n)\to c<\infty,\,\,J'(u_n)\to 0. $$ If $$
c<c^*=\frac{1}{N}S^{\frac{N}{2}}\min\bigg\{\Big(1-\frac{A+h(0)}{\L_N}\Big)^{\frac{N-1}{2}},\,\,
\Big(1-\frac{A+h(\infty)}{\L_N}\Big)^{\frac{N-1}{2}}\bigg\},$$
then $\{u_n\}$ has a converging subsequence.
\end{Theorem}
\begin{pf}
Let $\{u_n\}$ be a Palais-Smale sequence for $J$, then according
to $(h2)$, $\{u_n\}$ is bounded in $\dob(\ren )$. Then, up to
a subsequence, $ i)\,\, u_n\rightharpoonup u_0\mbox{ in } \dob
(\ren)$, $ii)\,\,  u_n\to u_0$ a.e., and $ iii)\,\,u_n\to u_0$ in
$L^\alpha _{loc}$, $\a \in [1,2^*)$. Therefore, by using the {\it
concentration compactness principle} by P. L. Lions, (see
\cite{PL1} and \cite{PL2}), there exists a subsequence (still denoted
by $\{u_n\}$) which
satisfies
\begin{enumerate}
\item $|\nabla u_n|^2\rightharpoonup d\mu\ge |\n u_0|^2+\sum_{j\in {\mathcal J} }\mu_j\delta_{x_j}+\mu _0\delta_0,$
\item $ |u_n|^{2^*}\rightharpoonup d\nu= |u_0|^{2^*}+\sum_{j\in {\mathcal J}}\nu_j\delta_{x_j}+\nu _0\delta_0,$
\item $  S\nu_j^{\frac2{2^*}}\le\mu_j \mbox{ for all } j\in{\mathcal J}\cup\{0\}$, where ${\mathcal J}$ is at most
countable,
\item $ \dfrac{u_n^2}{|x|^2}\rightharpoonup d \g= \dfrac{u_0^2}{|x|^2}+\g _0\delta_0,$
\item $\L _N\g_0\le\mu_0$.
\end{enumerate}
To study the  concentration at infinity of   the sequence we will
also need to consider the following  quantities $$
\nu_{\infty}=\dyle\lim_{R\to \infty}\limsup_{n\to
\infty}\int_{|x|>R}|u_n|^{2^*}dx,\quad \mu_{\infty}=\dyle\lim_{R\to
\infty}\limsup_{n\to \infty}\int_{|x|>R}|\n u_n|^{2^*}dx$$ and $$
\g_{\infty}=\dyle\lim_{R\to \infty}\limsup_{n\to
\infty}\int_{|x|>R}\frac{u_n^2}{|x|^2}dx.$$ We claim that ${\mathcal
J}$ is finite and for $j\in {\mathcal J}$,   either $\nu_j=0$ or
$\nu_j\ge S^{N/2}$. We follow closely the arguments in
\cite{AGP2}. Let $\e
>0$ and let $\phi $ be a smooth  cut-off function centered at $x_j$, $0\le \phi(x)\le
1$ such that
$$ \phi(x)= \left\{
\begin{array}{l}
1,\,\mbox{ if }|x-x_j|\le \dfrac{\e }{2}, \\ 0,\,\mbox{ if
}|x-x_j|\ge \e,
\end{array}
\right. $$ and $|\nabla \phi |\le \dfrac{4}{\e } $.
So we get
$$\begin{array}{lll} 0&=\lim\limits_{n\to \infty }\langle
J'(u_n),u_n\phi\rangle  \\ &  \\ &=\lim\limits_{n\to \infty
}\left(\dyle\irn |\nabla u_n|^2\phi+\irn u_n \nabla u_n \nabla
\phi - \irn\dfrac{A+ h(x)}{|x|^2}u_n^2 \phi -\irn \phi
|u_n|^{2^*}\right).
\end{array}
$$ From 1), 2) and 4) and since $0\notin \hbox{supp}( \phi )$ we
find that $$ \lim_{n\to\infty}\irn |\nabla u_n|^2\phi=\irn \phi d
\mu,\,\,\, \lim_{n\to\infty}\irn |u_n|^{2^*} \phi =\irn \phi d \nu
$$ and $$\lim\limits_{n\to\infty}\int_{B_{\e }(x_j)} \dfrac{A+
h(x)}{|x|^2}u_n^2 \phi=\int_{B_{\e }(x_j)} \dfrac{A+
h(x)}{|x|^2}u_0^2 \phi. $$ Taking limits as  $\e \to 0$ we obtain,
$$ \lim_{\e \to 0}\lim_{n \to\infty}\bigg|\irn u_n \nabla u_n
\nabla \phi\bigg|\to 0 .$$ Hence, $$ 0=\lim_{\e\to 0}\lim_{n\to
\infty }\langle J'(u_n),u_n\phi \rangle= \mu_j -\nu_j.$$ By 3) we
have that $ S\nu_j^{\frac2{2^*}}\le\mu_j$, then we obtain that
either $\nu_j=0$ or $\nu_j\ge S^{N/2}$, which implies that
${\mathcal J}$ is finite. The claim is proved.

\noindent Let us now study the possibility of concentration at $x=0$ and at $\infty$.
Let $\psi$ be a regular function such that $0\le \psi(x)\le 1$,
$$\psi(x)=
\left\{
\begin{array}{l}
1,\,\mbox{ if }|x|>R+1 \\ 0,\,\mbox{ if }|x|<R,
\end{array}
\right. $$ and $|\nabla \psi |\le \dfrac{4}{R} $. From \eqref{eq:minimiz} we obtain
that
\begin{equation}\label{eq:minimiz1}
\frac{\dyle\inr |\nab (u_n\psi)|^2dx -(A+h(\infty))\inr
\dfrac{\psi^2u^2_n}{|x|^2}dx}{\Big(\irn|\psi
u_n|^{2^*}\Big)^{2/2^*}}\ge
S\Big(1-\frac{A+h(\infty)}{\L_N}\Big)^{\frac{N-1}{N}}.
\end{equation}
Hence
$$
 \dyle\inr |\nab (u_n\psi)|^2dx -(A+h(\infty))\inr
\dfrac{\psi^2u^2_n}{|x|^2}dx\ge
S\Big(1-\frac{A+h(\infty)}{\L_N}\Big)^{\frac{N-1}{N}}\Big(\irn|\psi
u_n|^{2^*}\Big)^{2/2^*}.
$$
Therefore we conclude that
\begin{eqnarray*}
& \dyle\inr \psi ^2|\nab u_n|^2dx +\inr u_n^2|\nab \psi|^2dx+2\inr
u_n\psi \nab u_n\n \psi dx  \\ &\hskip2cm \ge  (A+h(\infty))\inr
\dfrac{\psi^2u^2_n}{|x|^2}dx+
S\Big(1-\dfrac{A+h(\infty)}{\L_N}\Big)^{\frac{N-1}{N}}\Big(\irn|\psi
u_n|^{2^*}\Big)^{2/2^*}.
\end{eqnarray*}
We claim that
$$\dyle\lim_{R\to \infty}\limsup_{n\to \infty}\bigg\{\inr
u_n^2|\nab \psi|^2dx+2\inr |u_n|\psi |\nab u_n||\n \psi|
dx\bigg\}=0.$$
Using H\"older inequality we obtain
\begin{eqnarray*}
\inr |u_n|\psi |\nab u_n||\n \psi| dx & \le &
\bigg(\int_{R<|x|<R+1} |u_n|^2|\n \psi|^2
dx\bigg)^{1/2}\bigg(\int_{R<|x|<R+1} |\n u_n|^2dx\bigg)^{1/2}.
\end{eqnarray*}
Hence \begin{eqnarray*} \limit\limits_{n\to \infty}\inr |u_n|\psi
|\nab u_n||\n \psi| dx & \le & C\bigg(\int_{R<|x|<R+1} |u_0|^2|\n
\psi|^2 dx\bigg)^{1/2}\\ &\le &  C \bigg(\int_{R<|x|<R+1}
|u_0|^{2^*}dx \bigg)^{2/2^*}\bigg(\int_{R<|x|<R+1} |\n \psi |^{N}dx
\bigg)^{2/N}\\ & \le & \overline{C} \bigg(\int_{R<|x|<R+1}
|u_0|^{2^*}dx\bigg)^{2/2^*}.
\end{eqnarray*}
Therefore we conclude that $$ \dyle\limit_{R\to
\infty}\limsup_{n\to \infty}\inr |u_n|\psi |\nab u_n||\n \psi|
dx\le \overline{C} \limit_{R\to \infty}\bigg(\int_{R<|x|<R+1}
|u_0|^{2^*}dx \bigg)^{2/2^*}=0. $$ Using the same argument we can
prove that $$ \dyle\lim_{R\to \infty}\limsup_{n\to \infty}\inr
u_n^2|\n \psi|^2=0. $$ Then we get
\begin{equation}\label{eq:infty}
\mu_{\infty}-(A+h(\infty))\g_{\infty}\ge
S\bigg(1-\frac{A+h(\infty)}{\L_N}\bigg)^{\frac{N-1}{N}}\nu_{\infty}^{2/2^*}.
\end{equation}
Since $\dyle\lim_{R\to \infty}\lim\limits_{n\to \infty }\langle
J'(u_n),u_n\psi\rangle=0$, we obtain that
$\mu_{\infty}-(A+h(\infty))\g_{\infty}\le\nu_{\infty}$. Therefore we
conclude that either $\nu_{\infty}=0$ or $\nu_{\infty}\ge
S^{\frac{N}{2}}\big(1-\frac{A+h(\infty)}{\L_N}\big)^{\frac{N-1}{2}}.$

\noindent The same holds for the concentration in $x_0=0$, namely that
either
$$
\nu_0=0 \mbox{  or
}\nu_0\ge
S^{\frac{N}{2}}\bigg(1-\frac{A+h(0)}{\L_N}\bigg)^{\frac{N-1}{2}}.
$$
As a conclusion we obtain
\begin{eqnarray*} c & = & J(u_n)-\frac{1}{2}\langle J'(u_n),u_n\rangle
  +o(1)\\ & = & \frac{1}{N}\irn
|u_n|^{2*}dx +o(1) = \frac{1}{N}\bigg\{\irn
|u_0|^{2*}dx+\nu_0+\nu_{\infty}+\sum_{j\in {\mathcal
J}}\nu_j\bigg\}.
\end{eqnarray*}
If we assume the existence of $j\in {\mathcal J}\cup\{0,\infty\}$ such that $\nu_j\neq 0$, then we obtain
that $c\ge c^*$ a contradiction with the hypothesis, then up to a subsequence $u_n\to u_0$ in
$\dob(\ren)$.
\end{pf}

\noindent  To find solutions  requires  to consider some path in
$\dob(\ren)$  along which the maximum of $J(\gamma(t))$ is less than
$c^*$. To do that, for $H=\max\{h(0),h(\infty)\}$,  we consider
$\{w_{\mu}\}$ the one parameter family of minimizer to problem
(\ref{eq:minimiz}) where $A$ is replaced by $A+H$. Then we have
the following result.
\begin{Theorem}\label{th:ex}
Suppose that $(h0)$, $(h1)$ and $(h2)$ hold. Assume the existence of
$\mu_0>0$ such that
\begin{equation}\label{eq:H}
\irn h(x)\frac{w^2_{\mu_0}(x)}{|x|^2}dx > H
\irn\frac{w^2_{\mu_0}(x)}{|x|^2}dx,\end{equation} then
(\ref{eq:poi1}) has at least a positive solution.
\end{Theorem}
\begin{pf}
Let $\mu_0$ be as in the hypothesis, then if we set
$$
f(t)=J(tw_{\mu_0})= \frac{t^2}{2}\bigg(\irn |\n
w_{\mu_0}|^2dx-\irn
\frac{A+h(x)}{|x|^2}w_{\mu_0}^2dx\bigg)-\frac{t^{2^*}}{2^*}\irn
|w_{\mu_0}|^{2^*}dx,\,t\ge 0
$$
we can see easily that $f$ achieves
its maximum at some $t_0>0$ and we can prove the existence of
$\rho>0$ such that $J(tw_{\mu_0})<0$ if $||tw_{\mu_0}||\ge \rho$.
By a simple calculation we obtain that
$$
t^{2^*-2}_0=\frac{\irn |\n
w_{\mu_0}|^2dx-\irn \frac{A+h(x)}{|x|^2}w^2_{\mu_0}dx}{\irn
|w_{\mu_0}|^{2^*}dx},
$$
and
\begin{equation}
J(t_0w_{\mu_0})=\max_{t\ge
0}J(tw_{\mu_0})=\frac{1}{N} \left(\frac{\irn |\n w_{\mu_0}|^2dx-\irn
\frac{A+h(x)}{|x|^2}w_{\mu_0}^2dx}{\big(\irn
|w_{\mu_0}|^{2^*}dx\big)^{2/2^*}}\right)^{N/2}.
\end{equation}
Using
(\ref{eq:H}) we obtain that
\begin{equation}
J(t_0w_{\mu_0})< \frac{1}{N} \left(\frac{\irn |\n
w_{\mu_0}|^2dx-(A+H)\irn \frac{w_{\mu_0}^2}{|x|^2}dx}{\big(\irn
|w_{\mu_0}|^{2^*}dx\big)^{2/2^*}}\right)^{N/2}=\frac{1}{N}S^{\frac{N}{2}}\bigg(1-\frac{A+H}{\L_N}\bigg)^{\frac{N-1}{2}}\le
c^* .
\end{equation}
We set $$ \G=\{\g\in C([0,1],\dob(\ren)):\,\,\g (0)=0\mbox{  and
}J(\g(1))<0\}. $$ Let $$ c=\dyle\inf_{\g\in \G }\max_{t\in
[0,1]}J(\g (t)). $$ Since $J(t_0w_{\mu_0})<c^*$, then  we get a
mountain pass critical point $u_0$. Then we have just to prove
that we can choose $u_0\ge 0$. We give two different proofs.

\noindent {\it First proof.} Consider the Nehari manifold,
\begin{align*}
M &\equiv  \{u\in \dob(\ren):\,\,u\neq 0 \mbox{  and }\langle
J'(u),u\rangle =0\}\\
 & =  \bigg\{u\in \dob(\ren):\,\,u\neq 0
\mbox{ and }\dyle\io |\n u|^2dx =\io \frac{A+h(x)}{|x|^2}u^2dx+\io
|u|^{2^*}dx\bigg\}.
\end{align*}
Notice that $u_0,\,\,|u_0|\in M$.
Since $u_0$ is a mountain pass solution to problem (\ref{eq:poi1})
then one can prove easily that $c\equiv J(u_0)=\dyle\min_{u\in
M}J(u)$ (see \cite{W}). Moreover as $J(|u_0|)=\dyle\min_{u\in
M}J(u)$, then $|u_0|$ is also a critical  point of $J$.

\noindent {\it Second proof.} Here we use a variation of the
deformation lemma. Since $u_0$ is a mountain pass critical point of $J$, which is even,
we have $$c=J(u_0)=J(|u_0|)=\max\limits_{t>0}J(t|u_0|).$$ Let
$t_1>0$ be such that $J(t_1|u_0|)<0$. We set
$\g_0(t)=t(t_1|u_0|)$ for $t\in [0,1]$. Notice that $\g_0\in \G$
and  $$ c=J(|u_0|)=\max_{t\in[0,1]}J(\g_0(t)).$$ If $|u_0|$ is a
critical point to $J$, then we have done. If not then using Lemma
3.7 of \cite{GO} we obtain that $\g_0$ can be deformed to a path
$\g_1\in \G$ with $\max_{t\in[0,1]}J(\g_1(t))<c$, a contradiction
with the definition of $c$ as a  {\it min-max value}.

\noindent Hence we have nonnegative solution to problem
\eqref{eq:poi1}. The positivity of the solution $u_0$ is an
application of the strong maximum principle by using
hypotheses $(h0)$ and $(h1)$.
\end{pf}

\noindent We give now some sufficient condition on $h$ to have
hypothesis (\ref{eq:H}).
\begin{Lemma}\label{eq:HH}
Suppose  one of the following hypotheses holds
\begin{enumerate}
\item[(1)] $h(x)\ge h(0)+c_1|x|^{\nu_{A+H}(N-2)}$ for $|x|$
small and $c_1>0$ if $h(0)\ge h(\infty)$, or
\item[(2)] $h(x)\ge h(\infty)+c_2|x|^{-\nu_{A+H}(N-2)}$ for $|x|$
large and $c_2>0$ if $h(\infty)\ge h(0)$,
\end{enumerate}
then there exists $\mu_0>0$ such that (\ref{eq:H}) holds.
\end{Lemma}
\begin{pf}Let $\d>0$ be small such that if $|x|<\d$ then
$h(x)\ge h(0)+c_1|x|^{\nu_{A+H}(N-2)}$. For simplicity of notation we
set $\nu_{A+H}=\nu$.
Let  $$I_{\d,\mu}=\dyle\int_{|x|<\d}(h(x)-H)\dfrac{dx}
{|x|^{(1-\nu)N+2\nu}(\mu^{2\nu}+|x|^{2\nu})^{N-2}},$$ then $$
I_{\d,\mu}\ge c_1\dyle\int_{|x|<\d}\dfrac{|x|^{\nu(N-2)} dx }
{|x|^{(1-\nu)N+2\nu}(\mu^{2\nu}+|x|^{2\nu})^{N-2}}.$$ Since $
\nu(N-2)-[(1-\nu)N+2\nu+2\nu(N-2)]=-N$, we conclude that
$I_{\d,\mu}\to \infty$ as $\mu\to 0$. On the other hand $$
\dyle\int_{|x|\ge \d}|h(x)-H|\dfrac{dx}
{|x|^{(1-\nu)N+2\nu}(\mu^{2\nu}+|x|^{2\nu})^{N-2}}\le
C\int_{|x|\ge \d}\dfrac{dx}{|x|^{(1+\nu)N-2\nu}}\le C(\d).$$
Therefore we get the existence of $\mu_0>0$ such that $$ \irn
(h(x)-H)\frac{w^2_{\mu_0}(x)}{|x|^2}dx\ge \int_{|x|<\d}
(h(x)-H)\frac{w^2_{\mu_0}(x)}{|x|^2}dx-\int_{|x|\ge \d}
|h(x)-H|\frac{w^2_{\mu_0}(x)}{|x|^2}dx>0.$$Then the result
follows.

\noindent The second case follows by using the same argument near
infinity.
\end{pf}
 \section{Perturbation of the nonlinear term: Multiplicity of positive solutions}
In this section we deal with the following problem
\begin{equation}\label{eq:poi11}
\begin{cases}
-\D u=\dfrac{\l}{|x|^2}u+k(x)u^{2^*-1}, \ x\in \ren,\\
u>0 \hbox{  in  }\ren, \mbox{  and  }u\in \dob(\ren),
\end{cases}
\end{equation}
where $N\ge 3$,\, $0<\l<\L_N$ and $k$ is a positive function.
\subsection{Existence}\label{ss:31}

Assume that $k$ verifies the hypothesis
$$
(K0) \quad k\in
L^{\infty}(\ren)\cap C(\ren)\,\,\hbox{ and } ||k||_{\infty}>\max
\{k(0), k(\infty)\} ,\hbox{ where }  k(\infty)\equiv
\limsup\limits_{|x|\to \infty}k(x).
$$
We associate to problem \eqref{eq:poi11}
the following functional
\begin{equation}
J_\l(u)= \frac12\irn |\n u|^2dx-\frac{\l}{2}\irn
\frac{u^2}{|x|^2}dx-\frac{1}{2^*}\irn k(x)|u|^{2^*}dx.
\end{equation}

\noindent As in the first section we have the following
Lemma.
\begin{Lemma}\label{lm:ps1}
Let $\{u_n\}\subset \dob(\ren)$ be a Palais-Smale sequence for
$J_{\l}$, namely $$ J_{\l}(u_n)\to c<\infty,\,\,J'_{\l}(u_n)\to 0.
$$
If
$$
c<\tilde c(\l)=\frac{1}{N}S^{\frac{N}{2}}\min\bigg\{||k||_{\infty}^{-\frac{N-2}{2}},
\,(k(0))^{-\frac{N-2}{2}}\Big(1-\frac{\l}{\L_N}\Big)^{\frac{N-1}{2}},
(k(\infty))^{-\frac{N-2}{2}}\Big(1-\frac{\l}{\L_N}\Big)^{\frac{N-1}{2}}\bigg\}
$$
then $\{u_n\}$ has a converging subsequence.
\end{Lemma}

\noindent The proof is similar to  the proof of  Theorem \ref{th:ps}.

\noindent In the case in which $k$ is a radial positive function, we
can prove the following improved Palais-Smale condition.
\begin{Lemma}\label{lm:psr}
Define
$$
\tilde c_1(\l)=\frac{1}{N}S^{\frac{N}{2}}\Big(1-\frac{\l}{\L_N}\Big)^{\frac{N-1}{2}}
\min\big\{(k(0))^{-\frac{N-2}{2}},\,(k(\infty))^{-\frac{N-2}{2}}\big\}.
$$
If
$\{u_n\}\subset \dob(\ren)$ is a Palais-Smale sequence for
$J_{\l}$,  namely
$$
J_{\l}(u_n)\to c, \quad J'_{\l}(u_n)\to 0,
$$
and $c<\tilde c_1$, then $\{u_n\}$ has a converging subsequence.
\end{Lemma}
\begin{remark}
This follows from the fact that the inclusion of $H^1_r(\O)\equiv
\{u\in L^2(\O):\,|\n u|\in L^2(\O), \,u \mbox{ radial}  \}$, where
$\O=\{x\in\ren:\,R_1<|x|<R_2\}$, in $L^q(\O)$ is compact for all
$1\le q<\infty$ and in particular for $q=2^*$, see \cite{K}.
\end{remark}
\noindent As a consequence we obtain the following existence
result.
\begin{Theorem}\label{th:ex1}
Let $k$ be a positive radial function such that
$(K0)$ is satisfied. Assume that there exists $\mu_0>0$ such that
\begin{equation}\label{eq:k}
\irn
k(x)w^{2^*}_{\mu_0}(x)dx
>\max\{k(0),k(\infty)\}
\irn w^{2^*}_{\mu_0}(x)dx,
\end{equation}
where $w_{\mu_0}$ is a
solution to problem
\begin{equation*}
\begin{cases}
-\D w=\dfrac{\l}{|x|^2}w+w^{2^*-1},\ x\in\ren,\\
w>0 \hbox{  in  }\ren, \mbox{  and  }w\in \dob(\ren).
\end{cases}
\end{equation*}
Then (\ref{eq:poi11}) has at least a positive solution.
\end{Theorem}
\begin{pf}
Since the proof is similar to the proof of Theorem \ref{th:ex}, we omit
it.
\end{pf}

\begin{remark}\label{r:3.5}
Assume that one of the following hypotheses holds
\begin{enumerate}
\item[(1)] $k(x)\ge k(0)+c_1|x|^{2\nu_{\l}}$ for $|x|$
small and $c_1>0$ if $k(0)\ge k(\infty)$, or
\item[(2)] $k(x)\ge k(\infty)+c_2|x|^{-2\nu_{\l}}$ for $|x|$
large and $c_2>0$ if $k(\infty)\ge k(0)$,
\end{enumerate}
then there exists $\mu_0>0$ such that (\ref{eq:k}) holds.
\end{remark}

\noindent Let us set
\[
b(\l)\equiv
\begin{cases}
+\infty &\text{if }k(0)=k(\infty)=0\\[5pt]
\min\bigg\{\displaystyle{(k(0))^{-\frac{N-2}{2}}\Big(1-\frac{\l}{\L_N}
\Big)^{\frac{N-1}{2}},(k(\infty))^{-\frac{N-2}{2}}\Big(1-\frac{\l}
{\L_N}\Big)^{\frac{N-1}{2}}}\bigg\}&\text{otherwise}.
\end{cases}
\]
\begin{Lemma}\label{l:eps0}
If $(K0)$ holds, there exists $\e_0>0$ such that
$||k||_{\infty}^{-\frac{N-2}{2}}\leq
b(\e_0)$ and
\begin{equation}\label{eq:ctilde}
\tilde c(\l)=\tilde c\equiv
\frac{1}{N}{S^{N/2}}{||k||_{\infty}^{-\frac{N-2}{2}}}
\end{equation}
for any
$0<\l\leq \e_0$.
\end{Lemma}
\begin{pf}
>From $(K0)$ it follows that if $\e_0$ is sufficiently small then
$||k||_{\infty}^{-\frac{N-2}{2}}\leq
b(\e_0)$ and hence from the definition of $\tilde c(\l)$ we obtain the result.
\end{pf}

\subsection{Multiplicity}\label{ss:32}

\noindent To find multiplicity results for problem
\eqref{eq:poi11} we need  the following extra hypotheses on $k$:
\begin{itemize}
\item[$(K1)$] the set ${\mathcal C}(k)=\Big\{a\in \ren\,\Big|\,
  k(a)=\max\limits_{x\in \ren}k(x)\Big\}$  is finite, say  ${\mathcal
    C}(k)=\big\{a_j\,|\,1\leq j\leq \text{\rm Card\,}({\mathcal C}(k))\big \}$;
\item[$(K2)$] there exists $2<\theta<N$ such that if $a_j\in {\mathcal C}(k)$ then
 $k(a_j)-k(x)=o(|x-a_j|)^{\theta}\mbox{ as }x\to a_j$.
\end{itemize}
\noindent Consider
$0<r_0\ll 1$ such that $B_{r_0}(a_j)\cap B_{r_0}(a_i)=\emptyset$ for $i
\neq j$, $1\leq i,j\leq \text{\rm Card\,}({\mathcal C}(k))$. Let
$\d=\frac{r_0}{3}$ and
for any $1\leq j\leq \text{\rm Card\,}({\mathcal C}(k))$ define the following
function
\begin{equation}\label{eq:tj}
T_j(u)=\dfrac{\irn\psi_j(x)|\n u|^2dx}{\irn |\n
u|^2dx}\,\, \hbox{  where  }\,\, \psi_j (x)=\min\{1,|x-a_j|\}.
\end{equation}
Notice that if $u\not\equiv 0$ and $T_j(u)\leq\delta$, then
\begin{eqnarray*}
r_0\int_{\ren\backslash B_{r_0}(a_j)}|\n u|^2dx & \le &
\int_{\ren\backslash B_{r_0}(a_j)}\psi_j(x)|\n u|^2dx\\ & \le &
\irn \psi_j(x)|\n u|^2dx\le \d \irn|\n u|^2dx=\frac{r_0}{3}\irn|\n
u|^2dx.
\end{eqnarray*}
Hence we have the following  property.
\begin{Lemma}\label{lm:estim1} Let $u\in\dob(\ren)$ be such that $T_j(u)\le \d$, then $$ \dyle\irn |\n u|^2\ge 3
\int_{\ren\backslash B_{r_0}(a_j)}|\n u|^2dx.$$
\end{Lemma}
\noindent As a consequence we obtain the following separation
result.
\begin{Corollary}\label{cor:separ} Consider $u\in \dob(\ren)$,  $u\not\equiv 0$,
such that $T_i(u)\le \d$ and $T_j(u)\le
\d$, then $i=j$.
\end{Corollary}
\begin{pf}
By  Lemma \ref{lm:estim1}  we obtain that
$$
2\irn|\n u|^2dx\ge 3\bigg(\int_{\ren\backslash B_{r_0}(a_i)}|\n
u|^2dx+\int_{\ren\backslash B_{r_0}(a_j)}|\n u|^2dx\bigg).
$$
If
$i\neq j$ we find that $$ 2\irn|\n u|^2dx\ge 3\irn|\n u|^2dx,$$
a contradiction if $u\not\equiv 0$.
\end{pf}

\noindent Consider the Nehari manifold,
\begin{equation}\label{eq:nehari}
M({\l})=\{u\in \dob(\ren):\,\,u\not\equiv 0 \mbox{
and }\langle J'_\l(u),u\rangle =0\}.
\end{equation}
Therefore if $u\in M({\l})$
$$\irn |\n u|^2dx-\l\irn
\frac{u^2}{|x|^2}dx=\irn k(x)|u|^{2^*}dx.
$$
Notice that for all
$u\in \dob(\ren)$,  $u\not\equiv 0$, there exists $t>0$
such that $tu\in M(\l)$ and for all $u\in M(\l)$ we have
$$
\irn
|\n u|^2dx-\l\irn \frac{u^2}{|x|^2}dx<(2^*-1)\irn k(x)|u|^{2^*}dx,
$$
hence, there exists $c_1>0$ such that
$$\forall u\in M(\l), \, \quad  ||u||_{\dob(\ren)}\ge c_1.$$
\begin{Definition}\label{def:level}
For any $0<\l<\Lambda_N$ and $1\leq j\leq \text{\rm Card\,}({\mathcal
  C}(k))$, let us consider
$$
M_j(\l)= \{u\in M(\l):\,\,T_j(u)<\d\}\mbox{ and its
boundary }\G_j(\l)=\{u\in M(\l):\,\,T_j(u)=\d\}.
$$
We define
$$
m_j(\l)=\inf\{J_\l(u):\,u\in M_j(\l)\} \mbox{ and
}\eta_j(\l)=\inf\{J_\l(u):\,u\in \G_j(\l)\}.
$$
\end{Definition}
The following two Lemmas give the behaviour of the functional with
respect to the critical level $\tilde c$.
\begin{Lemma}\label{lm:lem1}
Suppose that $(K0)$, $(K1)$, and $(K2)$ hold, then $M_j(\l)\neq\emptyset$ and
there exists $\e_1>0$ such that
\begin{equation}\label{eq:sss}
m_j(\l)<\tilde c\mbox{\quad  for all } 0<\l\le
\e_1\text{ and } 1\leq j\leq \text{\rm Card\,}({\mathcal C}(k)).
\end{equation}
\end{Lemma}
\begin{pf}
We set
\begin{equation}\label{eq:talenti}
v_{\mu,j}(x)=\frac{1}{(\mu^2+|x-a_j|^{2})^{\frac{N-2}{2}}}\mbox{\quad
and\quad  }u_{\mu,j}=\dfrac{v_{\mu,j}}{||v_{\mu,j}||_{{2^*}}},
\end{equation}
then $||u_{\mu,j}||_{{2^*}}=1$ and $\irn |\n u_{\mu,j}|^2dx =S$.
If
$$
t_{\mu,j}(\l)=\left(\dfrac{\irn |\n
u_{\mu,j}|^2dx-\l\irn|x|^{-2}{u^2_{\mu,j}}dx}{\irn
k(x)|u_{\mu,j}|^{2^*}dx}\right)^{\frac{N-2}{4}},
$$
 then
$t_{\mu,j}(\l)u_{\mu,j}\in M(\l)$. Making the change of variable
$x-a_j=\mu y$, we obtain
$$
T_j(t_{\mu,j}(\l)u_{\mu,j}) =
\dfrac{\irn\psi_j(x)|\n u_{\mu,j}|^2dx}{\irn |\n u_{\mu,j}|^2dx}=
\dfrac{\irn\psi_j(a_j+\mu y)|\n u_0(y)|^2dy}{\irn |\n
u_0(y)|^2dy}, $$
 where $u_0(x)$ is $u_{\mu,j}$ to scale
$\mu=1$ and concentrated in the origin.
Then $$\limit_{\mu\to 0}
T_j(t_{\mu,j}(\l)u_{\mu,j}) = \dfrac{\irn\psi_j(a_j)|\n
u_{0}(y)|^2dy}{\irn |\n u_{0}(y)|^2dy}=\psi_j(a_j)=0, $$ uniformly
in $\l$.
Hence we get the existence of $\mu_0$ independent of $\l$ such that if
$\mu <\mu_0$, then $t_{\mu,j}(\l)u_{\mu,j}\in M_j(\l)$. Notice that
\[
t_{\mu,j}(\l)\geq t_1(\l)\equiv \bigg(\|k\|_{\infty}^{-1}\Big(1-\frac{\l}{\Lambda_N}\Big)S\bigg)^{\frac{N-2}4}.
\]
In order to prove \eqref{eq:sss}, it is sufficient to
show the existence of $\mu<\mu_0$ such that if $0<\l<\e_1$ then
$$
\dyle\max_{t\geq t_1(\l)}J_\l(tu_{\mu,j})=J_\l(t_{\mu,j}(\l)u_{\mu,j})<\tilde c.$$
We have
$$
\max_{t\geq t_1(\l)}J_\l(tu_{\mu,j})\le
\max_{t>0}\Big\{\frac{t^2}{2}\irn |\n
u_{\mu,j}|^2dx-\frac{t^{2^*}}{2^*}\irn
k(x)|u_{\mu,j}|^{2^*}dx\Big\}-\frac 12 \l t_1^2(\l)\irn
\frac{u^2_{\mu,j}}{|x|^2}dx
$$
and
$$
\max_{t>0}\bigg\{\frac{t^2}{2}\irn |\n
u_{\mu,j}|^2dx-\frac{t^{2^*}}{2^*}\irn
k(x)|u_{\mu,j}|^{2^*}dx\bigg\}=\frac{1}{N}\bigg(\frac{S}{\big(\irn
k(x)|u_{\mu,j}|^{2^*}dx\big)^{2/{2^*}}}\bigg)^{N/2}.
$$
In view of assumption $(K2)$ we have that for some positive constants
$\bar c_1,\bar c_2$
\begin{align*}
\irn k(x)&|u_{\mu,j}|^{2^*}dx  =  ||k||_\infty-\irn
(k(a_j)-k(x))|u_{\mu,j}|^{2^*}dx \\
& =
||k||_\infty-\bar c_1\mu^{N}\irn\frac{k(a_j)-k(x)}{(\mu^2+|x-a_j|^2)^N}dx\\
& \ge
||k||_\infty-\bar c_1\mu^{N}\bigg\{\int_{B_{\d}(a_j)}\frac{\bar c_2|x-a_j|^{\theta}dx}
{(\mu^2+|x-a_j|^2)^N}+
2\|k\|_{\infty}\int_{\ren\backslash
  B_{\d}(a_j)}\frac{dx}{(\mu^2+|x-a_j|^2)^N}\bigg\}\\
& \ge
||k||_\infty-\bar c_1\mu^{N}\bigg\{\mu^{\theta-N}\bar c_2\int_{\ren}\frac{|y|^{\theta}dy}
{(1+|y|^2)^N}+
2\|k\|_{\infty}\int_{|y|\geq\delta}\frac{dy}{|y|^{2N}}\bigg\}\\
&= ||k||_\infty+O(\mu^{\theta}).
\end{align*}
Then we obtain that
\begin{eqnarray*}
\max_{t\geq t_1(\l)}J(tu_{\mu,j})& \le &
\frac{1}{N}\frac{S^{N/2}}{||k||_{\infty}^{\frac{N-2}{2}}+O(\mu^{\theta})}-
\frac12\l t_1^2(\l) \irn
\frac{u^2_{\mu,j}}{|x|^2}dx\\ & \le &
\frac{1}{N}\frac{S^{N/2}}{||k||_{\infty}^{\frac{N-2}{2}}}+O(\mu^{\theta})-\frac12 \l t_1^2(\l)\irn
\frac{u^2_{\mu,j}}{|x|^2}dx.
\end{eqnarray*}
Using estimate A.6 from \cite{SM} we obtain that for some positive
constant $c$
$$
\irn
\frac{u^2_{\mu,j}}{|x|^2}dx\ge c\mu^2\mbox{  as  }\mu\to 0.
$$
Therefore we get
\begin{eqnarray*}
\max_{t\geq t_1(\l)}J(tu_{\mu,j})& \le &
\frac{1}{N}\frac{S^{N/2}}{||k||_{\infty}^{\frac{N-2}{2}}}+O(\mu^{\theta})-
\frac12 c\l t_1^2(\l) \mu^2\\ & \le &
\tilde c+\bar c_3\mu^{\theta}- \frac12 c\l t_1^2(\l) \mu^2,
\end{eqnarray*}
where $\bar c_3$ is a positive constant. Since from $(K2)$ we have
$2<\theta<N$, we get the existence of $\e_1$ and $\mu_0$ such that
if $\mu<\mu_0$ and $0<\l<\e_1$, then $ \max\limits_{t\geq
t_1(\l)}J(tu_{\mu,j})<\tilde c$  and the result follows.
\end{pf}

\noindent We prove now the next result.
\begin{Lemma}\label{lm:estim3}
Suppose that $(K0)$, $(K1)$, and $(K2)$ are satisfied, then there exists $\e_2$ such
that for all $0<\l<\e_2$ we have $$ \tilde c<\eta_j(\l).$$
\end{Lemma}
\begin{pf} We argue by contradiction. We assume the existence of $\l_n\to 0$ and
$\{u_n\}$ such that $u_n\in\G_j(\l_n)$ and $J_{\l_n}(u_n)\to c\le
\tilde c=\frac{1}{N}{S^{N/2}}{||k||_\infty^{-\frac{N-2}{2}}}$. We can
easily prove that $\{u_n\}$ is bounded. Then up to a subsequence
we get the existence of $l>0$ such that $$ \limit\limits_{n\to
\infty}\irn |\n u_n|^2dx=\limit\limits_{n\to \infty}\irn
k(x)|u_n|^{2^*}dx=l.$$ Notice that $l\ge
{S^{N/2}}{||k||_\infty^{-\frac{N-2}{2}}}$. On the other hand, by the
definition of $\{u_n\}$ we have,
\begin{eqnarray*}
\frac{1}{N}l+o(1)& = & J_{\l_n}(u_n)\\ & = & \frac{1}{2}\irn |\n
u_n|^2dx-\frac{\l_n}{2}\irn \frac{u^2_n}{|x|^2}-\frac{1}{2^*}\irn
k(x)|u_n|^{2^*}dx \\ & \le &
\frac{1}{N}{S^{N/2}}{||k||_\infty^{-\frac{N-2}{2}}}+o(1).
\end{eqnarray*}
Then we conclude that $l={S^{N/2}}{||k||_\infty^{-\frac{N-2}{2}}}$
and therefore,
\begin{equation}\label{eq:main}
\limit\limits_{n\to \infty}\irn
\big(||k||_{{\infty}}-k(x)\big)|u_n|^{2^*}dx=0.
\end{equation}
We set $w_n=\dfrac{u_n}{||u_n||_{ {2^*}}}$, then
$||w_n||_{{2^*}}=1$ and $$\limit_{n\to \infty}\irn|\n
w_n|^2dx =S.$$
Hence by using the concentration compactness arguments by P.L. Lions
(see also Proposition 5.1 and 5.2 in \cite{S}), we get the existence
of $w_0\in
\dob(\ren)$ such that $w_n$ converges to $w_0$ weakly in $\dob(\ren)$
(up to a subsequence) and one of the following alternatives holds
\begin{enumerate}
\item $w_0\not\equiv 0$ and   $w_n\to w_0$ strongly in the $\dob(\ren)$.
\item  $w_0\equiv 0$ and either
\begin{itemize}
\item[i)]  $|\nabla w_n|^2\rightharpoonup d\mu=S\d_{x_0}$ and $|w_n|^{2^*}\to
d\nu=\d_{x_0}$

or
\item[ii)]  $|\nabla w_n|^2\rightharpoonup d\mu_{\infty}=S\d_{\infty}$ and $|w_n|^{2^*}\rightharpoonup
d\nu_{\infty}=\d_{\infty}$.
\end{itemize}
\end{enumerate}
The last case means that
$$
\nu_{\infty}=\dyle\lim_{R\to
\infty}\limsup_{n\to \infty}\int_{|x|>R}|w_n|^{2^*}dx=1\mbox{ and
} \mu_{\infty}=\dyle\lim_{R\to \infty}\limsup_{n\to
\infty}\int_{|x|>R}|\n w_n|^{2^*}dx=S.
$$
If the first alternative
holds,  from (\ref{eq:main}) we obtain that $$
\limit\limits_{n\to \infty}\irn
\big(||k||_{{\infty}}-k(x)\big)w^{2^*}_ndx= \irn
\big(||k||_{{\infty}}-k(x)\big)w^{2^*}_0dx=0,
$$
a contradiction with the
fact that $k$ is not a constant.

\noindent Assume that we have the alternative 2 i), then since
$T_j(w_n)=T_j(u_n)=\d$, we conclude that $$
\d=T_j(w_n)=\limit_{n\to \infty}\dfrac{\irn\psi_j(x)|\n
w_n|^2dx}{\irn |\n w_n|^2dx}=\psi_j(x_0).$$ Hence the
concentration is impossible in any point $a_j\in{\mathcal
  C}(k)$. On the other hand from
(\ref{eq:main}) we obtain that
$$
0=\limit\limits_{n\to
\infty}\irn
\big(||k||_{{\infty}}-k(x)\big)w^{2^*}_ndx=||k||_\infty-k(x_0),
$$
a contradiction.

\noindent To analyze concentration at $\infty$, consider a regular
function $\xi$ satisfying
$$
\xi(x)= \left\{
\begin{array}{l}
1,\,\mbox{ if }|x|>R+1 \\ 0,\,\mbox{ if }|x|<R,
\end{array}
\right.
$$
where $R$ is chosen in a such way that $|a_j|<R-1$ for
all $j$. Then we have
\begin{eqnarray*}
\d=T_j(w_n)& = &
\limit_{n\to \infty}\dfrac{\irn\psi_j(x)|\n w_n|^2dx}{\irn |\n
w_n|^2dx}\\ & = & \limit_{R\to \infty}\limit_{n\to
\infty}\dfrac{\irn\xi (x)|\n w_n|^2dx +\irn (1-\xi(x))\psi_j(x)|\n
w_n|^2dx}{\irn |\n w_n|^2dx}.
\end{eqnarray*}
Since $\limit_{n\to \infty}\irn (1-\xi(x))\psi_j(x)|\n
w_n|^2dx=0$, we conclude that $$\d=\limit_{R\to
\infty}\limit_{n\to \infty}\dfrac{\irn\xi (x)|\n w_n|^2dx}{\irn
|\n w_n|^2dx}=1,$$ a contradiction if we choose $\d<1$.
So we conclude.
\end{pf}

\noindent We need now the following Lemma that is suggested by the work of Tarantello \cite{TA}.
See also \cite{CC1}.
\begin{Lemma}\label{lm:ta} Assume that $0<\l<\min\{\e_1,\e_2\}$ where $\e_1,\,\e_2$ are
given by  Lemmas \ref{lm:lem1} and \ref{lm:estim3}. Then for
all $u\in M_j(\l)$ there exists $\rho_u>0$ and a differentiable
function $$f:B(0,\rho_u)\subset \dob(\ren)\to \re$$
 such that $f(0)=1$ and for all $w\in
B(0,\rho_u)$ we have $f(w)(u-w)\in M_j(\l)$. Moreover for all
$v\in \dob(\ren)$ we have
\begin{equation}\label{main2}
\langle f'(0),v\rangle=-\dfrac{\displaystyle{2\irn \n u\n vdx -2\l\irn \dfrac{u
v}{|x|^2}dx -2^*\irn k(x)|u|^{2^*-2}uvdx}}{\displaystyle{\irn |\n u|^2dx
-\l\irn\dfrac{u^2}{|x|^2}dx -(2^*-1)\irn k(x)|u|^{2^*}dx}}.
\end{equation}
\end{Lemma}
\begin{pf}
Let $u\in M_j(\l)$ and let $G: \re \times \dob(\ren)\to \re$ be
the function defined by
$$
G(t,w)=t\bigg(\irn |\n (u-w)|^2dx
-\l\irn\frac{(u-w)^2}{|x|^2}dx\bigg) -t^{2^*-1}\irn
k(x)|u-w|^{2^*}dx.
$$
Then $G(1,0)=0$ and $G_t(1,0)=\irn |\n u|^2dx
-\l\irn\frac{u^2}{|x|^2}dx -(2^*-1)\irn k(x)|u|^{2^*}dx\neq 0$
(since $u\in M_j(\l)$). Then by using the Implicit Function
Theorem we get the existence of $\rho_u>0$ small enough and of a
differentiable function $f$ satisfying the required property.
Moreover, notice that
$$
\langle f'(0),v\rangle=-\dfrac{\langle
G_w(1,0),v\rangle}{G_t(1,0)}=-\dfrac{\displaystyle{2\irn \n u\n vdx -2\l\irn
\dfrac{u v}{|x|^2}dx -2^*\irn k(x)|u|^{2^*-2}uvdx}}{\displaystyle{\irn |\n u|^2dx
-\l\irn\dfrac{u^2}{|x|^2}dx -(2^*-1)\irn k(x)|u|^{2^*}dx}}.
$$
\end{pf}

\noindent We are now in position to prove the main result.
\begin{Theorem}\label{th:mutlip}
Assume that $(K0)$, $(K1)$, and $(K2)$ hold, then there exists $\e_3$ small such that
for all $0<\l<\e_3$ equation (\ref{eq:poi11}) has $\text{\rm
  Card\,}({\mathcal C}(k))$ positive
solutions $u_{j,\l}$ such that
\begin{equation}\label{eq:main3} |\n u_{j,\l}|^2\to
{S^{N/2}}{||k||_\infty^{-(N-2)/2}}\d_{a_j}\mbox{  and }|
u_{j,\l}|^{2^*}\to {S^{N/2}}{||k||_\infty^{-N/2}}\d_{a_j}\mbox{ as
}\l\to 0.
\end{equation}
\end{Theorem}
\begin{pf}
Assume that $0<\l<\e_3=\min\{\e_0,\e_1,\e_2\}$, where $\e_0$, $\e_1$
and $\e_2$ are
given by  Lemmas \ref{l:eps0}, \ref{lm:lem1} and \ref{lm:estim3}. Let $\{u_n\}$ be a
minimizing sequence to $J_{\l}$ in $M_j(\l)$,  that is, $u_n\in M_j(\l)$ and
$J_{\l}(u_n)\to m_j(\l)$ as $n\to \infty$. Since
$J_{\l}(u_n)=J_{\l}(|u_n|)$, we can choose $u_n\ge 0$. Notice
that we can prove the existence of $c_1, c_2$ such that $c_1\le
||u_n||_{\dob(\ren)}\le c_2$. By the Ekeland variational principle
we get the existence of a subsequence denoted also by $\{u_n\}$
such that $$ J_{\l}(u_n)\le m_j(\l)+\frac{1}{n}\mbox{  and  }
J_{\l}(w)\ge J_\l(u_n)-\frac{1}{n}||w-u_n||\mbox{  for all }w\in
M_j(\l).  $$ Let $0<\rho<\rho_{n}\equiv \rho_{u_n}$ and $f_n\equiv f_{u_n}$, where
$\rho_{u_n}$ and $f_{u_n}$ are given by Lemma \ref{lm:ta}. We set $v_\rho=\rho v$
where $||v||_{\dob(\ren)}=1$, then $v_\rho\in B(0,\rho_n)$ and we
can apply Lemma \ref{lm:ta} to obtain that
$w_\rho=f_n(v_\rho)(u_n-v_\rho)\in M_j(\l)$.

\noindent Therefore we get
\begin{eqnarray*}
\frac{1}{n}||w_\rho-u_n|| & \ge & J_{\l}(u_n)-J_\l(w_\rho)=\langle
J'_{\l}(u_n), u_n-w_\rho\rangle +o(||u_n-w_\rho||)\\ & \ge & \rho
f_n(\rho v)\langle J'_{\l}(u_n), v\rangle +o(||u_n-w_\rho||).
\end{eqnarray*}
Hence we conclude that $$ \langle J'_{\l}(u_n), v\rangle\le
\frac{1}{n}\frac{||w_\rho-u_n|| }{\rho f_n(\rho v)}(1+o(1))\hbox{ as  } \rho\to 0.$$
Since $ |f_n(\rho v)|\to |f_n(0)|\ge c $ as $\rho\to 0$ and
\begin{eqnarray*}
\frac{||w_\rho-u_n|| }{\rho} & = & \frac{||f_n(0)u_n-f_n(\rho
v)(u_n-\rho v)||}{\rho }\\ & \le &
\frac{||u_n|||f_n(0)-f_n(\rho v)|+|\rho||f_n(\rho v)|}{\rho}\le
C|f'_n(0)|||v||+c_3\le c.
\end{eqnarray*}
Then we conclude that $J'_\l(u_n)\to 0$ as $n\to \infty$.
Hence $\{u_n\}$ is a Palais-Smale sequence for $J_\l$. Since
$m_j(\l)<\tilde c$ and $\tilde c=\tilde c(\l)$ for $\l\leq \e_0$, then
from Lemma \ref{lm:ps1} we get the existence result.

\noindent To prove (\ref{eq:main3}) we follow the proof of Lemma \ref{lm:estim3}. Assume $\l_n\to 0$ as
$n\to \infty$ and let $u_n\equiv u_{j_0,\l_n}\in M_{j_0}(\l_n)$ be a solution to problem (\ref{eq:poi11})
with $\l=\l_n$. Then up to a subsequence we get the existence of $l_1>0$ such that $$ \limit\limits_{n\to
\infty}\irn |\n u_n|^2dx=\limit\limits_{n\to \infty}\irn k(x)|u_n|^{2^*}dx=l_1.$$ Therefore as in the
proof of Lemma \ref{lm:estim3} we obtain that $l_1=S^{N/2}||k||_\infty^{-\frac{N-2}{2}}$ and $$
\limit\limits_{n\to \infty}\irn (||k||_{\infty}-k(x))u^{2^*}_ndx=0. $$ We set
$w_n=\dfrac{u_n}{||u_n||_{2^*}}$, then $||w_n||_{2^*}=1$ and $\dyle\limit_{n\to \infty}||
w_n||_{\dob(\ren)}^2=S$. Hence we get the existence of $w_0\in \dob(\ren)$ such that one of the following
alternatives holds
\begin{enumerate}
\item $w_0\not\equiv 0$ and   $w_n\to w_0$ strongly in the $\dob(\ren)$.
\item  $w_0\equiv 0$ and either
\begin{itemize}
\item[i)]  $|\nabla w_n|^2\rightharpoonup d\mu=S\d_{x_0}$ and $|w_n|^{2^*}\to
d\nu=\d_{x_0}$

or
\item[ii)]  $|\nabla w_n|^2\rightharpoonup d\mu_{\infty}=S\d_{\infty}$ and $|w_n|^{2^*}\rightharpoonup
d\nu_{\infty}=\d_{\infty}$.
\end{itemize}
\end{enumerate}
As in Lemma \ref{lm:estim3}, the  alternative 1  and the
alternative 2 ii) do not hold. Then we conclude that the unique possible behaviour  is
the  alternative 2. i), namely, we get the existence of $x_0\in
\ren$ such that $$ |\nabla w_n|^2\rightharpoonup
d\mu=S\d_{x_0}\mbox{  and  } |w_n|^{2^*}\rightharpoonup
d\nu=\d_{x_0}.$$ Since
\begin{eqnarray*} \irn |\n w_n|^2dx & = & S+o(1)=S\irn
|w_n|^{2^*}dx +o(1)=\frac{S}{||k||_{\infty}}\irn
k(x)|w_n|^{2^*}dx+o(1)\\ & = &
\frac{S}{||k||_{\infty}}k(x_0)+o(1),\end{eqnarray*} then we obtain
that $x_0\in {\mathcal C}(k)$. Using Corollary \ref{cor:separ}, we
conclude that $x_0=a_{j_0}$ and the result follows.
\end{pf}
\begin{remark}
As in \cite{CC3}, we can prove the same kind of results under more general condition on $k$.
For instance, we
can assume that $k$ changes sign and the following conditions hold,
\begin{itemize}
\item[$(K'1)$] $\max\limits_{x\in \ren}k(x)>0$ and ${\mathcal C'}(k)=\{a\in \ren\,|\, k(a)=\max\limits_{x\in \ren}k(x)\}$
is a finite set.
\item[$(K'2)$] $(K2)$ holds.
\end{itemize}
In this case the level at which the Palais-Smale conditions fails becomes
$$
\hat c(\l)=\frac{1}{N}S^{\frac{N}{2}}\min\bigg\{||k_+||_{\infty}^{-\frac{N-2}{2}},
\,(k_+(0))^{-\frac{N-2}{2}}\Big(1-\frac{\l}{\L_N}\Big)^{\frac{N-1}{2}},
(k_+(\infty))^{-\frac{N-2}{2}}\Big(1-\frac{\l}{\L_N}\Big)^{\frac{N-1}{2}}\bigg\}.
$$
\end{remark}
\section{Category setting.}
In this section we use the Lusternik-Schnirelman category theory
to get multiplicity results for problem (\ref{eq:poi11}), we refer to
\cite{Am} for a complete discussion. We follow the argument by
Musina  in \cite{Mu}. We assume that $k$ is a nonnegative function and
that $0<\l<\overline{\e}_0$ where $\overline{\e}_0$ is chosen in a
such way that
$\big(1-\frac{\overline{\e}_0}{\L_N}\big)^{\frac{N-1}2}>\frac12$
and $\overline{\e}_0\le \e_0$, being $\e_0$ given in Lemma
\ref{l:eps0}. We set for $\delta>0$
$$
{\mathcal C}(k)=\{a\in \ren\,|\,
k(a)=||k(x)||_\infty\}\mbox{ and } {\mathcal C}_{\d}(k)=\{x\in
\ren: dist(x,{\mathcal C}(k))\le \d\}.
$$
We suppose
that $(K2)$ holds and
$$ (K3)\quad \mbox{there exist
}R_0,\,d_0>0\mbox{  such that }\sup\limits_{|x|>R_0}|k(x)|\le
||k||_{\infty}-d_0.$$ Let
$M(\l)$ be defined by \eqref{eq:nehari}. Consider
$$\tilde{M}(\l)\equiv \{u\in
M(\l):\,J_\l(u)<\tilde c\}.$$ Then we have the following local
Palais-Smale condition.
\begin{Lemma}\label{lm:psi}
Let $\{v_n\}\subset M(\l)$ be such that
\begin{equation}\label{eq:portatil} J_\l(v_n)\to c<\tilde{c}\mbox{
and  }{J'_\l}_{|_{M(\l)}}(v_n)\to 0,
\end{equation}
then $\{v_n\}$ contains a
converging subsequence.
\end{Lemma}
\begin{pf}
Assume that $\{v_n\}$ satisfies \eqref{eq:portatil}, then
there exists $\{\a_n\}\subset \re$ such that
\begin{equation}\label{eq:estim3}
J'_\l(v_n)-\a_n G'_\l(v_n)\to 0 \mbox{  as }n\to \infty \mbox{  in
}\mathcal{D}^{-1,2}(\ren)
\end{equation}
where $G_{\l}(u)=\langle J'_\l(u), u\rangle$.
Since $\{v_n\}\subset M(\l)$ and
$J_\l(v_n)\le \tilde c$, we have
$r_1\le ||v_n||_{\dob(\ren)}\le r_2$ for some constants $r_1,r_2>0$.
Using $v_n$ as a test function in
\eqref{eq:estim3} we conclude that $\a_n\to 0$ as $n\to \infty$.
Hence $\{v_n\}$ is a Palais-Smale sequence for $J_\l$ at the level
$c<\tilde c$ and then the result follows by using Lemma
\ref{lm:ps1}.
\end{pf}

\noindent To prove that $\tilde{M}(\l)\neq\emptyset$ we give the next
result.
\begin{Lemma}\label{lm:cat1}
There exists $\overline{\e}_1>0$ such that if
$0<\l<\l_0:=\min\{\overline{\e}_0,\overline{\e}_1\}$, then $\tilde{M}(\l)\neq
\emptyset$. Moreover for any $\{\l_n\}\subset\re_+$ such that
$\l_n\to 0$ as
$n\to \infty$ and $\{v_n\}\subset
\tilde{M}(\l_n)$, there exist $\{x_n\}\subset \ren$ and $\{r_n\}\subset \re_+$ such that $x_n\to x_0\in {\mathcal C}(k)$,\,$r_n\to 0$ as
$n\to \infty$ and
\begin{equation}\label{eq:cat2}
v_n-\bigg(\frac{S}{||k||_\infty}\bigg)^{\frac{N-2}{4}}u_{r_n}
(\cdot-x_n)\to
0\mbox{ in }\dob(\ren),
\end{equation}
where
\begin{equation}\label{talenti}
u_{r}(x)=\dfrac{C_r}{(r^2+|x|^{2})^{\frac{N-2}{2}}}\end{equation}
and $C_r$ is the normalizing constant to be $||u_r||_{2^*}=1$.
\end{Lemma}
\begin{pf}
The first assertion follows by using the same argument as in Lemma
\ref{lm:lem1} since we have $$ \max_{t>0}J_{\l}(tw_{\mu,x}) \le
\frac{1}{N}\frac{S^{N/2}}{||k||_{\infty}^{\frac{N-2}{2}}}+O(\mu^{\theta})-
c\l\mu^2<\tilde c \ \mbox{ for }\mu \mbox{ small and }2<\theta <N,$$
where $w_{\mu,x}(y)=\dfrac{C}{(\mu^2+|y-x|^{2})^{\frac{N-2}{2}}}$,
$x\in \mathcal{C}(k)$ and $C$ is the normalizing constant such
that $||w_{\mu,x}||_{2^*}=1$ (see also \cite{BN}). As a
consequence, there exists $\l_0$ such that for all $0<\l<\l_0$ the
set $\tilde{M}(\l)$ is not empty. To prove the second part of the
Lemma, eventually passing to a subsequence we set $$
\limit\limits_{n\to \infty}\irn |\n v_n|^2dx=\limit\limits_{n\to
\infty}\irn k(x)|v_n|^{2^*}dx=l.$$ Then as in Lemma
\ref{lm:estim3} we can prove that
$l={S^{N/2}}{||k||_\infty^{-\frac{N-2}{2}}}$ and
\begin{equation}\label{eq:categ1}
\limit\limits_{n\to \infty}\irn
(||k||_{{\infty}}-k(x))v^{2^*}_ndx=0.
\end{equation}
Consider the normalized function $w_n=\dfrac{v_n}{||v_n||_{ {2^*}}}$
and $$\limit_{n\to \infty}\irn|\n w_n|^2dx
=S.$$
  Using the concentration-compactness arguments  by P.L. Lions,
we obtain the existence of $\{x_n\}\subset \ren$ and $\{r_n\}\subset \re_+$ such that
\begin{equation}\label{eq:cat02}
w_n-u_{r_n}(\cdot-x_n)\to 0\mbox{ in }\dob(\ren),
\end{equation}
and $w_n\db w_0\in \dob(\ren)$. Moreover by the same argument as in the proof of Lemma \ref{lm:estim3}
the weak limit is $w_0=0$. We will show now that  the concentration
 at infinity is not possible. Indeed if concentration at $\infty$ occurs, by using (\ref{eq:categ1}) and $(K3)$ we obtain
\begin{eqnarray*}
||k||_{\infty}& = & \irn
k(x)w^{2^*}_ndx+o(1)=\int_{\ren\backslash
B_{R_0}(0)}k(x)w^{2^*}_ndx+o(1)\\ & \le &
\sup\limits_{|x|>R_0}|k(x)|\int_{\ren\backslash
B_{R_0}(0)}w^{2^*}_ndx+o(1)\le (||k||_{\infty}-d_0)+o(1),
\end{eqnarray*}
which is a contradiction.  Then the unique possible
concentration  is at some point $x_0\in\ren$. Hence  we conclude
that, up to a subsequence,
$r_n\to 0$ and $$ |\n u_{r_n}(x-x_n)|^2\db S\d_{x_0}.$$ Using
(\ref{eq:categ1}) it is easy to obtain that $x_0\in  {\mathcal C}(k)$.
\end{pf}
\begin{remark} Notice that as a consequence of the above Lemmas we
obtain the existence of at least $cat(\tilde{M}(\l))$ solutions
that eventually can change sign.
\end{remark}

\noindent
Hereafter we concentrate our study on the analysis  of
$cat(\tilde{M}(\l))$, the behaviour of the energy,  and
the positivity of solutions.

\noindent If $R_0$ is like in hypothesis $(K3)$, we define
\[
\xi(x)=
\begin{cases}
x    &\mbox{ if }|x|\le  R_0,\\
R_0\dfrac{x}{|x|}& \mbox{ if
}|x|\ge R_0,
\end{cases}
\]
and for $u\in \dob(\ren)$ such that $u\neq 0$ we set
\begin{equation}\label{eq:beta}
\Xi(u)=\frac{\irn \xi(x)|\n u|^2dx}{\irn |\n u|^2dx}.
\end{equation}
We recall that for $u\in \dob(\ren)$ such that $u\neq 0$ we have
$t_{\l}(u)u\in M(\l)$ where $t_{\l}(u)$ is given by
$$t_{\l}(u)=\Bigg(\dfrac{\irn|\n
u|^2dx-\l\irn\frac{u^2}{|x|^2}dx}{\irn
k(x)|u|^{2^*}dx}\Bigg)^{\frac{N-2}{4}}.$$ Let $\Psi_\l:\ren\to
\dob(\ren)$ be given by $$ \Psi_\l
(x)=t_{\l}(u_{\mu_\l}(\cdot-x))u_{\mu_\l}(\cdot-x),$$ where
$u_{\mu_\l}$ is given by \eqref{eq:talenti}, $\mu_{\l}\equiv
g(\l)$ such that $g(\l)\to 0$ as $\l \to 0$. Notice that if $x\in
{\mathcal{C}}(k)$ and $\l$ is sufficiently small, then
\begin{equation}\label{eq:autonoma} J_\l(\Psi_{\l} (x))=
\max_{t>0}J_{\l}(tu_{\mu_\l}(\cdot-x))
\le
\frac{1}{N}\frac{S^{N/2}}{||k||_{\infty}^{\frac{N-2}{2}}}+O(\mu_\l^{\theta})-
c\l\mu^2_\l<\tilde c. \end{equation}
Then  we can prove the
existence of $\l_0, c_1, c_2>0$ such that for all $0<\l<\l_0$ we have
$\Psi_\l(x)\in \tilde{M}(\l)$, $J_\l(\Psi_{\l}(x))=\tilde c+o(1)$ as $\l\to0$, and
$c_1<t_{\l}(u_{\mu_\l}(\cdot-x))<c_2$ for all $x\in {\mathcal C}(k)$. As a
consequence, taking limits   for $\l\to 0$ we obtain by Lemma
\ref{lm:cat1} that for any $x\in {\mathcal{C}}(k)$
\begin{equation}\label{eq:catestim}
|\n \Psi_\l(x)|^2\db
d\mu=S^{N/2}||k||_\infty^{-\frac{N-2}{2}}\d_{x}\mbox{  and   } |
\Psi(x)|^{2^*}\db d\nu=(S||k||_\infty^{-1})^{N/2}\d_{x}.
\end{equation}
We prove now the next result.
\begin{Lemma}\label{lm:antes}
For $\l\to 0$ we have
\begin{enumerate}
\item $\Xi(\Psi_\l(x))=x+o(1)$ uniformly for $x\in B_{R_0}(0)$.
\item  $\sup\{dist(\Xi(u), \mathcal{C}(k)): u\in \tilde{M}(\l)\}\to 0.$
\end{enumerate}
\end{Lemma}
\begin{pf}
Let $x\in B_{R_0}(0)$, then by \eqref{eq:catestim} we obtain that
$$
\Xi(\Psi_\l(x))=\frac{\irn \xi(y)|\n \Psi_{\l}(x)|^2dy}{\irn |\n
\Psi_\l(x) |^2dy}=\frac{\irn \xi(y)d\mu}{\irn d\mu}+o(1)=x+o(1)\quad\text{as }\l\to0.
$$
To
prove the second assertion we take $\l_n\to 0$ and let $v_n\in
\tilde{M}(\l_n)$, then by Lemma \ref{lm:cat1} we get the existence of
$\{x_n\}\subset \ren$ and $\{r_n\}\subset \re_+$ such that such
that $x_n\to x_0\in{\mathcal C}(k)$,\,$r_n\to 0$ as $n\to \infty$
and
$$
v_n-\bigg(\frac{S}{||k||_\infty}\bigg)^{\frac{N-2}{4}}u_{r_n}(\cdot-x_n)\to
0\mbox{ in }\dob(\ren).
$$
Since $\Xi$ is a continuous function we
obtain that $$ \Xi(v_n)=\frac{\irn \xi(x)|\n v_n|^2dx}{\irn |\n
v_n|^2dx}=\frac{\irn \xi(x)|\n u_{r_n}(\cdot-x_n)|^2dx}{\irn |\n
u_{r_n}(\cdot-x_n)|^2dx}+o(1)=\xi(x_0)+o(1).$$ Since $x_0\in
{\mathcal C}(k)\subset B_{R_0}(0)$ we conclude that $\xi(x_0)=x_0$
and the result follows.
\end{pf}

\noindent We are now able to prove the main result.
\begin{Theorem}\label{th:ultimo}
Assume that hypotheses $(K0)$, $(K2)$ and $(K3)$ hold and let
$\d>0$. Then there exists $\l_0>0$ such that for all
$0<\l<\l_0$, equation \eqref{eq:poi11} has at least $cat_{{\mathcal
C}_{\d}(k)}{\mathcal C}(k)$ solutions.
\end{Theorem}
\begin{pf}
Given $\d>0$ there exists
$\l_0(\d)>0$ such that by  Lemma \ref{lm:antes} and \eqref{eq:autonoma},
 for $0<\l<\l_0(\d)$ we have that $\Psi_\l(x)\in
\tilde M(\l)$ for any $x\in {\mathcal C}(k)$, and
 $$
 |\Xi(\Psi_\l(x))-x|<\d\mbox{  for all }x\in
B_{R_0}(0)\mbox{ and }\Xi(u)\in {\mathcal C}_{\d}(k)\mbox{  for
all }u\in \tilde{M}(\l). $$ Let ${\mathcal
H}(t,x)=x+t(\Xi(\Psi_\l(x))-x)$ where $(t,x)\in [0,1]\times
{\mathcal C}(k)$, then ${\mathcal H}$ is a continuous function and
$dist({\mathcal H}(t,x),{\mathcal C}(k))\le \d$ for all $(t,x)\in
[0,1]\times {\mathcal C}(k)$. Hence ${\mathcal H}([0,1]\times
{\mathcal C}(k))\subset {\mathcal C}_{\d}(k)$.

\noindent Since ${\mathcal H}(0,x)=x$ and ${\mathcal
H}(1,x)=\Xi(\Psi_\l(x))$, then we conclude that $\Xi\circ\Psi_\l$ is
homotopic to the inclusion ${\mathcal C}(k)\hookrightarrow {\mathcal C}_\d(k)$.
Since $J_\l$ satisfies the Palais-Smale condition below the level
$\tilde c$, to prove the Theorem we need just to prove that
$cat (\tilde{M}(\l))\ge cat_{{\mathcal C}_{\d}(k)}{\mathcal C}(k)$.

\noindent Suppose that $\{M_i\}$, $i=1,...,n_0$, is a closed covering
of $\tilde{M}(\l)$, then for any $i=1,...,n_0$ there exists a homotopy $${\mathcal
H_i}:[0,1]\times M_i\to \tilde{M}(\l)$$ such that $$ {\mathcal
H_i}(0,u)=u\mbox{  for all }u\in M_i\mbox{  and  }{\mathcal
H_i}(1,\cdot)=\mbox{constant for }i=1,...,n_0.$$ Notice that from
\eqref{eq:autonoma}, we obtain that $\Psi_\l({\mathcal
C}(k))\subset \tilde{M}(\l)$. We set ${\mathcal
C_i}=\Psi_{\l}^{-1}(M_i)$, then ${\mathcal C_i}$ is closed in
${\mathcal C}_{\delta}(k)$ and ${\mathcal C}(k)\subset \cup_i{\mathcal
C_i}\subset {\mathcal C}_\delta(k)$. Then we have just only to show that ${\mathcal C_i}$ are
contractible in ${\mathcal C}_{\d}(k)$. We set ${\mathcal
G_i}:[0,1]\times {\mathcal C_i}\to {\mathcal C}_{\d}(k)$ where
${\mathcal G_i}(t,x)=\Xi({\mathcal H_i}(t,\Psi_{\l}(x)))$. Then
$$
{\mathcal G_i}(0,x)=\Xi\circ \Psi_{\l}(x)\mbox{  for all } x\in {\mathcal C_i}\mbox{
 and  }{\mathcal G_i}(1,\cdot)=\mbox{ constant for }i=1,...,n_0.$$
Since $\Xi\circ\Psi_\l$ is
homotopic to the inclusion ${\mathcal C}(k)\hookrightarrow {\mathcal C}_\d(k)$ we have that ${\mathcal C_i}$ are contractible in ${\mathcal C}_{\d}(k)$. To complete the proof it remains  to prove
that any solution has a fixed sign. We follow the argument used in
\cite{CL}. Assume that $u=u^+-u^-$ with $u^+\ge 0$,\,$u^-\ge 0$
and $u^+\not\equiv 0$,\,$u^-\not\equiv 0$. Then we have
\begin{equation}\label{eq:positv1}
\begin{array}{lll}
& \dyle\irn |\n u^\pm|^2dx -\l\irn \frac{|u^\pm|^2}{|x|^2}dx \ge
S(1-\frac{\l}{\L_N})^{\frac{N-1}{N}}\Big(\irn |u^\pm|^{2^*}dx
\Big)^{2/{2^*}}\\& \ge \dyle
S(1-\frac{\l}{\L_N})^{\frac{N-1}{N}}||k||^{-\frac{2}{2^*}}_\infty\Big(\irn
k(x)|u^\pm|^{2^*}dx\Big )^{2/{2^*}}
\end{array}
\end{equation}
Since $u$ is a solution to problem \eqref{eq:poi11} we obtain that
\begin{equation}\label{eq:positv2}
\irn |\n u^\pm|^2dx -\l\irn \frac{|u^\pm|^2}{|x|^2}dx = \irn
k(x)|u^\pm|^{2^*}dx.
\end{equation}
Therefore we conclude that
\begin{eqnarray*}
\tilde c & > & J_\l(u)=\frac{1}{N}\irn k(x)|u|^{2^*}dx
=\frac{1}{N}\bigg\{\irn k(x)|u^+|^{2^*}dx+\irn k(x)|u^-|^{2^*}dx\bigg\}\\ &
\ge &
\frac{2S}{N}^{\frac{N}{2}}\bigg(1-\frac{\l}{\L_N}\bigg)^{\frac{N-1}{2}}
||k||_\infty^{-\frac{N-2}{2}}.
\end{eqnarray*}
Hence we obtain $2(1-\frac{\l}{\L_N})^{\frac{N-1}{2}}\le 1$ which contradicts
 the choice of $\l$.
\end{pf}
\begin{remark}
\
\begin{itemize}
\item[i)]If ${\mathcal{C}}(k)$ is finite, then for $\l$ small,
  equation \eqref{eq:poi11} has at
least $\text{Card}(\mathcal{C}(k))$ solutions.
\item[ii)] We give now a typical example where equation \eqref{eq:poi11} has
infinity many solutions. Let $\eta: \re \to \re_+$ such that $\eta$
is regular, $\eta (0)=0$ and $\eta(r)=1$ for $r\ge \frac12$. We
define $k_1$ on $[0,1]\subset \re$ by
$$
k_1(r)=
\begin{cases}
0 &\mbox{  if  }r=0,\\
\eta(r)(1-|\sin(\frac{1}{r-\frac12})|^{\theta})&\mbox{ if }0<r\le
1,\\
\end{cases}
$$
where $2<\theta<N$.  Notice that $k_1$ has infinitely
many global maximums archived on the set
$$\mathcal{C}(k_1)=\bigg\{r_n=\frac12+\frac{1}{n\pi}\mbox{ for }n\ge
1\bigg\}.$$ Now we define $k$ to be any continuous bounded function
such that $k(x)=k_1(|x|)$ if $|x|\le 1$, $||k||_\infty \le 1$ and
$\limit\limits_{|x|\to \infty}k(x)=0$. Since for all $m\in \ene $,
there exists $\d(m)$ such that
$\text{cat}_{\mathcal{C}_\d(k)}(\mathcal{C}(k))=m$, then we conclude
that equation \eqref{eq:poi11} has at least $m$ solutions for
$0<\l<\l(\d)$.
\item[iii)] Let us note that if $\delta$ becomes larger, then
  $\text{cat}_{\mathcal{C}_\d(k)}(\mathcal{C}(k))$ decreases, so that
  Theorem \ref{th:ultimo} is interesting for $\delta$ small.
\end{itemize}

\end{remark}



\end{document}